\documentclass[12pt]{article}
\usepackage{amsmath}
\usepackage{amsthm}
\usepackage{amsfonts}
\usepackage{epic,eepic}
\let\oldlabel=\label
\def\prellabel{\marginparsep=1em\marginparwidth=44pt
    \def\label##1{\oldlabel{##1}\ifmmode\else\ifinner\else
         \marginpar{{\footnotesize\ \\ \tt
                    ##1}}\fi\fi}}
\def\NZQ{\mathbb}               
\def\NN{{\NZQ N }}

\def\ZZ{{\NZQ Z}}

\def\opn#1#2{\def#1{\operatorname{#2}}} 
\opn\chara{char}
\opn\gr{gr}
\opn\Rees{{\mathcal R}}
\opn\rank{rank}
\def\pot#1#2{#1[\kern-0.28ex[#2]\kern-0.28ex]}

\let\sect=\cap
\let\dirsum=\oplus

\let\Sect=\bigcap
\let\Dirsum=\bigoplus

\let\to=\rightarrow

\def\Implies{\ifmmode\Longrightarrow \else
     \unskip${}\Longrightarrow{}$\ignorespaces\fi}
\def\implies{\ifmmode\Rightarrow \else
     \unskip${}\Rightarrow{}$\ignorespaces\fi}

\let\:=\colon
\let\epsilon=\varepsilon
\let\phi=\varphi
\let\kappa=\varkappa
\def\email#1{{\tt #1}}
\def\SetSize{\fontsize{12}{14.4}\selectfont}

\makeatletter

\renewcommand\section{\@startsection {section}{1}{\z@}%
                                   {-3.5ex \@plus -1ex \@minus -.2ex}%
                                   {2.3ex \@plus.2ex}%
                                   {} }
\renewenvironment{thebibliography}[1]
     {\vspace{3.5ex \@plus 1ex \@minus .2ex}
      \noindent REFERENCES\par\vspace{2.3ex \@plus.2ex}%
      \list{\@biblabel{\@arabic\c@enumiv}}%
           {\settowidth\labelwidth{\@biblabel{#1}}%
            \leftmargin=\labelwidth
            \@openbib@code
            \usecounter{enumiv}%
            \let\p@enumiv\@empty
            \renewcommand\theenumiv{\@arabic\c@enumiv}}%
      \sloppy
      \clubpenalty4000
      \@clubpenalty \clubpenalty
      \widowpenalty4000%
      \sfcode`\.\@m}
     {\def\@noitemerr
       {\@latex@warning{Empty `thebibliography' environment}}%
      \endlist}
\renewcommand{\@biblabel}[1]{#1.}
\makeatother

\newtheoremstyle{theorem}
  {}
  {}
  {\itshape}
  {}
  {}
  {.}
  {.5em}
  {}

\newtheoremstyle{definition}
  {}
  {}
  {}
  {}
  {\MakeUppercase}
  {.}
  {.5em}
  {}

\theoremstyle{theorem}
\newtheorem{theorem}{THEOREM}[section]
\newtheorem{lemma}[theorem]{LEMMA}
\newtheorem{corollary}[theorem]{COROLLARY}
\newtheorem{proposition}[theorem]{PROPOSITION}

\theoremstyle{definition}
\newtheorem{remark}[theorem]{REMARK}

\newtheorem{examples}[theorem]{EXAMPLES}
\newtheorem{definition}[theorem]{DEFINITION}
\newtheorem{question}[theorem]{QUESTION}
\newtheorem{questions}[theorem]{QUESTIONS}

\opn\ini{in}
\opn\KRS{KRS}

\opn\krs{krs}
\opn\diag{diag}
\opn\DD{{\mathcal D}}
\opn\SS{{\mathcal S}}
\opn\MM{{\mathcal M}}
\opn\GL{GL}
\def\w{w}
\opn\height{height}
\opn\length{length}
\def\sep{\,|\,}
\def\addots{\mathinner{\mkern1mu\raise1pt\hbox{.}\mkern2mu\raise4pt\hbox{.}
        \mkern2mu\raise7pt\vbox{\kern7pt\hbox{.}}\mkern1mu}}

\def\Box#1#2#3{\multiput(#1,#2)(1,0){2}{\line(0,1)1}
                           \multiput(#1,#2)(0,1){2}{\line(1,0)1}
                           \put(#1,#2){\makebox(1,1){$#3$}}}
\def\EBox#1#2{\Box#1#2{}}
\def\DBox#1#2#3{\multiput(#1,#2)(1,0){2}{\dashline{0.2}(0,0)(0,1)}
                           \multiput(#1,#2)(0,1){2}{\dashline{0.2}(0,0)(1,0)}
                           \put(#1,#2){\makebox(1,1){$#3$}}}
\def\EDBox#1#2{\DBox#1#2{}}
\def\LBox#1#2#3#4{\multiput(#1,#2)(#4,0){2}{\line(0,1)1}
                           \multiput(#1,#2)(0,1){2}{\line(1,0){#4}}
                           \put(#1,#2){\makebox(#4,1){$#3$}}}

\def\Ci#1#2{\put(#1.5,#2.5){\circle{0.7}}}
\def\Cross#1#2{\put(#1.5,#2.5){\line(-1,-1){0.5}}\put(#1.5,#2.5){\line(-1,1){0.5
}}

\put(#1.5,#2.5){\line(1,-1){0.5}}\put(#1.5,#2.5){\line(1,1){0.5}}}

\begin{document}
\thispagestyle{empty} \vspace*{1.5in} {\fontsize{14}{16.8}\selectfont \noindent
\uppercase{KRS and determinantal ideals\par}}

\SetSize
\vspace{2\baselineskip}

\noindent\uppercase{Winfried Bruns}, Universit\"at Osnabr\"uck, FB
Mathematik/Informatik, 49069 Osna\-br\"uck, Germany,
\email{Winfried.Bruns@mathematik.uni-osnabrueck.de}\\[1\baselineskip]
\uppercase{Aldo Conca}, Dipartimento di Matematica e Fisica, Universit\'a di
Sassari, Via Vienna 2, 07100 Sassari, Italy, \email{conca@ssmain.uniss.it}

\vspace{4\baselineskip plus 1 \baselineskip minus 1\baselineskip}

\section{\uppercase{Introduction}}\label{Intro}

Let $K$ be a field and $X$ an $m\times n$ matrix of indeterminates. The
determinantal ideals in $K[X]$ are the ideals $I_t$ generated by the $t$-minors
of $X$, $1\le t\le \min(m,n)$, and ideals related to them.

The Knuth--Robinson--Schensted correspondence (KRS) is a powerful tool for the
computation of Gr\"obner bases of determinantal ideals. For this purpose it has
first been used by Sturmfels \cite{Stu}. Then Herzog and Trung \cite{HT} have
considerably extended the class of ideals to which KRS can be applied. In a
different direction Sturmfels' method has been generalized by Bruns and Conca
\cite{BC} and Bruns and Kwieci\'nski \cite{BK}. While Herzog and Trung use
Gr\"obner bases in order to derive numerical results, the papers \cite{BC} and
\cite{BK} aim at structural information, mainly on powers of determinantal
ideals and the corresponding Rees algebras.

The crucial point in the application of KRS to Gr\"obner bases  is to show the
equality $\ini(I)=\KRS(I)$ for the ideals $I$ under consideration. We call
these ideals \emph{in-KRS}. Here $\ini(I)$ is the initial ideal of $I$ with
respect to a so-called diagonal term order on $K[X]$, and $\KRS(I)$ is the
image of $I$ under the automorphism of the polynomial ring $K[X]$ induced by
$\KRS$ -- in the strict sense KRS is a bijection from the set of standard
bitableaux (or standard monomials) $\SS$ to the set of monomials $\MM$ of
$K[X]$. Both $\SS$ and $\MM$ are $K$-bases of $K[X]$: for $\SS$ this is
asserted by the straightening law of Doubilet--Rota--Stein. Since $\ini(I)$ is
a monomial ideal, one must assume that $I$ has a basis of standard bitableaux.

The first three sections of the paper are an expanded version of the first
author's lecture in the conference. Section \ref{SectStr} recapitulates the
straightening law, and Section \ref{SectKRS} introduces KRS. Section
\ref{SectGr} explains the common ideas underlying the results of \cite{BC} and
\cite{BK}. For this purpose we develop a conceptual framework in which KRS
\emph{invariants} play the central role. Such invariant is a function $F:\DD\to
\NN$ defined on the set $\DD$ of all bitableaux (or products of minors) that,
roughly speaking, is compatible with the straightening law and, moreover,
satisfy the condition
$$
F(\Sigma)=\max\{F(\Delta): \Delta\in\DD,\ \ini(\Delta)=\KRS(\Sigma)\}.
$$
It is then easy to see that each of the ideals $I_k(F)$ generated by all
standard bitableaux $\Sigma$ with $F(\Sigma)\ge k$ satisfies the condition
$\ini(I_k(F))=\KRS(I_k(F))$. Even more is true: $I_k(F)$ is G-KRS, i.~e.\ in
addition to being in-KRS, $I_k(F)$ has a Gr\"obner basis of bitableaux. The
class of G-KRS ideals is closed under sums, that of in-KRS ideals $I$ is closed
under sums and intersections, and therefore one then obtains many G-KRS or at
least in-KRS ideals.

It has been shown in \cite{BC} that the functions $\gamma_t$ introduced by De
Concini, Eisenbud and Procesi \cite{DEP} are KRS invariants. This fact allows
one to compute the Gr\"obner bases, or at least the initial ideals of the
symbolic powers of the $I_t$ and products $I_{t_1}\cdots I_{t_s}$. The family
$\alpha_k$ of KRS invariants found by Greene has been used in \cite{BK} for the
analysis of the ideal underlying MacPherson's graph construction in the generic
case.

In Section \ref{SectSh} we show that all ideals that are generated by products
of minors and do not ``prefer any rows or columns'' of the matrix $X$ are
in-KRS, at least if $\chara K$ is $0$ or $>\min(m,n)$. In characteristic $0$
this is exactly the class of ideals that have a standard monomial basis and are
stable under the natural action of $\GL(m,K)\times \GL(n,K)$ on $K[X]$. In fact
all these ideals can be written as sums of intersections of symbolic powers of
the ideals $I_t$, and the symbolic powers are G-KRS, as stated above.

Section \ref{SectBas} characterizes those among all the ideals of Section
\ref{SectSh} that are even G-KRS. We show that these are essentially the sums
of the ideals $J(k,d)$ introduced in \cite{BK} and for which Greene's theorem
yields the property of being G-KRS. Since each KRS invariant can be derived
from a family of G-KRS ideals, this shows that Greene's functions $\alpha_k$
are truly basic KRS invariants, at least if one considers functions
$F:\DD\to\NN$ for which $F(\Delta)$ only depends on the shape of $\Delta$.

Section \ref{SectIn} complements the results of \cite{BC}. We show that the
formation of initial ideal and symbolic power commute for the ideals $I_t$.
This result can be interpreted as a description of the semigroup of monomials
in the initial algebra of the symbolic Rees algebra by linear inequalities.

In Section \ref{SectCo}  we turn to a potential new KRS invariant
$\gamma_\delta$ related with the ideal $I(X,\delta)$ cogenerated by a minor
$\delta$. Except in the case in which $I_\delta=I_t$, these do not only depend
on shape and therefore constitute an interesting new class of functions. Though
\cite{HT} gives some information on $\gamma_\delta$, we have not yet been able
to show that these are KRS invariants.

The application of KRS to determinantal ideals has also been investigated by
Abhyankar and Kulkarni \cite{AK1,AK2}. Furthermore, variants of the KRS can be
used to study ideals of symmetric matrices  of indeterminates (Conca \cite{Co})
or ideals generated by Pfaffians of alternating matrices (\cite{HT}, Bae\c tica
\cite{Ba}, De Negri \cite{DeN}).

There are now excellent discussions of KRS available in textbooks; see Fulton
\cite{F} and Stanley \cite{Sta}.

\section{\uppercase{The straightening law}}\label{SectStr}

Let $K$ be a field and $X$ an $m\times n$ matrix of indeterminates over $K$.
For a given positive integer $t\leq \min(m,n)$, we consider the ideal
$I_t=I_t(X)$ generated by the $t$-minors (i.~e.\ the determinants of the
$t\times t$ submatrices) of $X$ in the polynomial ring $R=K[X]$ generated by
all the indeterminates $X_{ij}$.

 From the viewpoint of algebraic geometry $R$ should be regarded as the
coordinate ring of the variety of $K$-linear maps $f\: K^m\to K^n$. Then
$V(I_t)$ is just the variety of all $f$ such that $\rank f<t$, and $R/I_t$ is
its coordinate ring.

The study of the determinantal ideals $I_t$ and the objects related to them has
numerous connections with invariant theory, representation theory, and
combinatorics. For a detailed account we refer the reader to Bruns and Vetter
\cite{BV}.

Almost all of the approaches one can choose for the investigation of
determinantal rings use standard bitableaux and the straightening law. The
principle governing this approach is to consider all the minors of $X$ (and not
just the $1$-minors $X_{ij}$) as generators of the $K$-algebra $R$ so that
products of minors appear as ``monomials''. The price to be paid, of course, is
that one has to choose a proper subset of all these ``monomials'' as a linearly
independent $K$-basis: the standard bitableaux are a natural choice for such a
basis, and the straightening law tells us how to express an arbitrary product
of minors as a $K$-linear combination of the basis elements. (In \cite{BC},
\cite{BK} and \cite{BV} standard bitableaux were called \emph{standard
monomials}; however, we will have to consider the ordinary monomials in $K[X]$
so often that we reserve the term monomial for products of the $X_{ij}$.)

In the following
$$
[a_1,\dots,a_t\sep b_1,\dots,b_t]
$$
stands for the determinant of the submatrix $(X_{a_ib_j}\: i=1,\dots,t,\
j=1,\dots,t)$.

The letter $\Delta$ always denotes a product $\delta_1\cdots\delta_w$ of
minors, and we assume that the sizes $|\delta_i|$ (i.~e.\ the number of rows of
the submatrix $X'$ of $X$ such that $\delta_i=\det(X')$) are \emph{descending},
$|\delta_1|\ge \dots\ge |\delta_w|$. By convention, the empty minor $[\sep]$
denotes $1$. The \emph{shape} $|\Delta|$ of $\Delta$ is the sequence
$(|\delta_1|,\dots,|\delta_w|)$. If necessary we may add factors $[\sep]$ at
the right hand side of the products, and extend the shape accordingly.

A product of minors is also called a \emph{bitableau}. The choice of this term
bitableau is motivated by the graphical description of a product $\Delta$ as a
pair of Young tableaux as in Figure \ref{Young}:
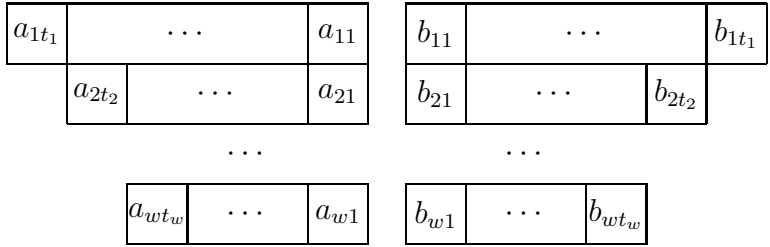
\begin{figure}[hbt]
\begin{gather*} \unitlength=0.8cm
\begin{picture}(6,4)(0,0)
\Box03{a_{1t_1}} \LBox13{\cdots}4 \Box53{a_{11}} \Box12{a_{2t_2}}
\LBox22{\cdots}3 \Box52{a_{21}} \put(2,1){\makebox(4,1){$\cdots$}}
\Box20{a_{wt_w}} \LBox30{\cdots}2 \Box50{a_{w1}}
\end{picture} \hspace{0.5cm}
\begin{picture}(6,4)(0,0)
\Box03{b_{11}} \LBox13{\cdots}4 \Box53{b_{1t_1}} \Box02{b_{21}}
\LBox12{\cdots}3 \Box42{b_{2t_2}} \put(0,1){\makebox(4,1){$\cdots$}}
\Box00{b_{w1}} \LBox10{\cdots}2 \Box30{b_{wt_w}}
\end{picture}
\end{gather*}
\caption{A bitableau}\label{Young}
\end{figure}
Every product of minors is represented by a bitableau and, conversely, every
bitableau stands for a product of minors if the length of the rows is
decreasing from top to bottom, the entries in each row are strictly increasing
from the middle to the outmost box, the entries of the left tableau are in
$\{1,\dots,m\}$ and those of the right tableau are in $\{1,\dots,n\}$. These
conditions are always assumed to hold.

For formal correctness one should consider the bitableaux as purely
combinatorial objects and distinguish them from the ring-theoretic objects
represented by them, but since there is no real danger of confusion, we simply
identify them.

Whether $\Delta$ is a standard bitableau is controlled by a partial order of
the minors, namely
\begin{multline*}
[a_1,\dots,a_t\sep b_1,\dots,b_t] \preceq [c_1,\dots,c_u\sep d_1,\dots,d_u]\\
\iff\quad t\ge u\quad\text{and}\quad a_i\le c_i,\ b_i\le d_i,\ i=1,\dots,u.
\end{multline*}
A product $\Delta=\delta_1\cdots\delta_w$ is called a \emph{standard bitableau}
if
$$
\delta_1\preceq\dots\preceq\delta_w,
$$
in other words, if in each column of the bitableau the indices are
non-decreasing from top to bottom. The letter $\Sigma$ is reserved for standard
bitableaux.

The fundamental straightening law of Doubilet--Rota--Stein says that every
element of $R$ has a unique presentation as a $K$-linear combination of
standard bita\-bleaux (for example, see Bruns and Vetter \cite{BV}):

\begin{theorem}\label{straight}
\begin{itemize}
\item[(a)] The standard bitableaux are a $K$-vector space basis of $K[X]$.
\item[(b)] If the product $\delta_1\delta_2$ of minors is not a standard
bitableau, then it has a representation
$$
\delta_1\delta_2=\sum x_i\epsilon_i\eta_i,\qquad x_i\in K,\ x_i\neq0,
$$
where $\epsilon_i\eta_i$ is a standard bitableau,
$\epsilon_i\prec\delta_1,\delta_2\prec\eta_i$ (here we must allow that
$\eta_i=1$).
\item[(c)] The standard representation of an arbitrary bitableau $\Delta$,
i.e.\ its representation as a linear combination of standard bitableaux
$\Sigma$, can be found by successive application of the straightening relations
in \emph{(b)}.
\item[(d)]  Moreover, at least one $\Sigma$ with $|\Sigma|=|\Delta|$ appears
with a non-zero coefficient in the standard representation of $\Delta$.
\end{itemize}
\end{theorem}

Let $e_1,\dots,e_m$ and $f_1,\dots,f_n$ denote the canonical $\ZZ$-bases of
$\ZZ^m$ and $\ZZ^n$ respectively.  Clearly $K[X]$ is a
$\ZZ^m\dirsum\ZZ^n$-graded algebra if we give $X_{ij}$ the `` vector bidegree''
$e_i\dirsum f_j$. All minors are homogeneous with respect to this grading, and
therefore the straightening relations must preserve the multiplicities with
which row and column indices occur on the left hand side.

The straightening law implies that the ideals $I_t$ have a $K$-basis of
standard bitableaux:

\begin{corollary}\label{Itstan}
The standard bitableaux $\Sigma=\delta_1\cdots\delta_w$ such that
$|\delta_1|\ge t$ form a $K$-basis of $I_t$.
\end{corollary}

All these standard bitableaux are elements of $I_t$ since $\delta_1\in I_t$ if
$|\delta_1|\ge t$. Conversely, every $x\in I_t$ can be written as a $K$-linear
combination of products $\delta M$ where $\delta$ is a minor of size $t$ and
$M$ is a monomial. Properties (b) and (c) of the straightening law imply
that the
standard bitableaux in the standard presentation of $\delta M$ have the
required property.

We say that an ideal $I\subset R$ has a \emph{standard basis} if $I$ is the
$K$-vector space spanned by the standard bitableaux $\Sigma\in I$.

\section{\uppercase{The Knuth--Robinson--Schensted correspondence}}\label{SectKRS}

Let $\Sigma$ be a standard bitableau. The Knuth--Robinson--Schensted
correspondence (see Fulton \cite{F} or Stanley \cite{Sta}) sets up a bijective
correspondence between standard bitableaux and monomials in the ring $K[X]$. We
use the version of KRS given by Herzog and Trung \cite{HT}.

If one starts from bitableaux, the correspondence is constructed from the
\emph{deletion} algorithm. Let $\Sigma=(a_{ij} | b_{ij} )$ be a non-empty
standard bitableau. Then
one constructs a pair of integers $(\ell,r)$ and a standard bitableau $\Sigma'$
as follows.
\begin{itemize}
\item[(a)]
One chooses the largest entry $\ell$ in the \emph{left} tableau of $\Sigma$;
suppose that $\{ (i_1,j_1),\allowbreak \dots,\allowbreak (i_u,j_u)\}$,
$i_1<\dots<i_u$, is the set of indices $(i,j)$ such that $\ell=a_{ij}$.
\item[(b)] Then the boxes at the pivot position $(p,q)=(i_u,j_u)$ in the
left and the
right tableau are removed.
\item[(c)] The entry $\ell=a_{pq}$ of the removed box in the left tableau is
the first component of the output, and $b_{pq}$ is stored in $s$.
 \item[(d)] The second and third components of the output are
determined by a ``push out'' procedure on the \emph{right} tableau as follows:
\begin{itemize}
\item[(i)]  if $p=1$, then $r=s$ is the second component of the output, and
the third is the standard bitableau $\Sigma'$ that has now been created;
\item[(ii)] otherwise the entry $b_{pq}$ is moved one row up and pushes out
the right most entry $b_{p-1k}$ such that $b_{p-1k}\le b_{pq}$ whereas
$b_{p-1k}$ is stored in $s$.
\item[(iii)] one replaces $p$ by $p-1$ and goes to step (i).
\end{itemize}
\end{itemize}
It is now possible to define KRS recursively: One sets $\KRS([\sep])=1$, and
$\KRS(\Sigma)=\KRS(\Sigma')X_{\ell r}$ for $\Sigma\neq[\sep]$.

We give an example in Figure \ref{Del}. The circles in the left tableau mark
the pivot position, those in the right mark the chains of ``push-outs'':
\begin{figure}[hbt]
$$
\begin{gathered}
\begin{picture}(4,2)(0,0)
\Box311 \Box213 \Box114 \Box015 \Box302 \Box206 \Ci20
\end{picture}
\hspace{0.5cm}
\begin{picture}(4,2)(0,0)
\Box011 \Box112 \Box213 \Box316 \Box004 \Box105 \Ci10 \Ci21
\end{picture}
\\
\begin{picture}(4,2)(0,0)
\Box311 \Box213 \Box114 \Box015 \Box302 \Ci01
\end{picture}
\hspace{0.5cm}
\begin{picture}(4,2)(0,0)
\Box011 \Box112 \Box215 \Box316 \Box004 \Ci31
\end{picture}
\\
\begin{picture}(4,2)(0,0)
\Box311 \Box213 \Box114 \Box302 \Ci11
\end{picture}
\hspace{0.5cm}
\begin{picture}(4,2)(0,0)
\Box011 \Box112 \Box215 \Box004 \Ci21
\end{picture}
\end{gathered}
\qquad\qquad\qquad
\begin{gathered}
\begin{picture}(2,2)(0,0)
\Box111 \Box013 \Box102 \Ci01
\end{picture}
\hspace{0.5cm}
\begin{picture}(2,2)(0,0)
\Box011 \Box112 \Box004 \Ci11
\end{picture}
\\
\begin{picture}(2,2)(0,0)
\Box111 \Box102 \Ci10
\end{picture}
\hspace{0.5cm}
\begin{picture}(2,2)(0,0)
\Box011 \Box004 \Ci00 \Ci01
\end{picture}
\\
\begin{picture}(2,2)(0,0)
\Box111 \Ci11
\end{picture}
\hspace{0.5cm}
\begin{picture}(2,2)(0,0)
\Box014  \Ci01
\end{picture}
\end{gathered}
$$
\caption{The KRS algorithm}\label{Del}
\end{figure}
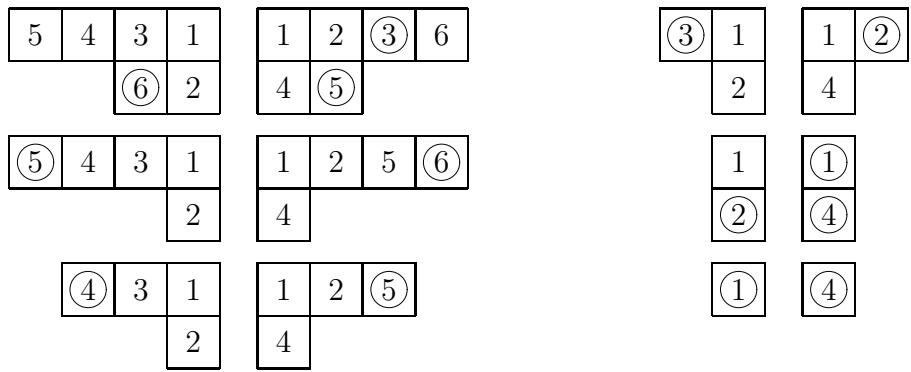
In this example we have
$$
\KRS(\Sigma)=X_{14}X_{21}X_{32}X_{45}X_{56}X_{}X_{63}.
$$
It is often more convenient to denote the output by a two row array instead of
a monomial by making the row indices of the factors the upper row and the
column indices the lower row; in the example
$$
\krs(\Sigma)=\begin{pmatrix}1&2&3&4&5&6\\4&1&2&5&6&3\end{pmatrix}.
$$
In general we set
$$
\krs(\Sigma)=\begin{pmatrix}u_1&\dots&u_w\\v_1&\dots&v_w\end{pmatrix}.
$$
In both rows indices may appear several times; however, the indices $u_i$ in
the upper row are non-decreasing from left to right, and if $u_i=u_{i+1}$, then
$v_i\ge v_{i+1}$ for the indices in the lower row, as is easily checked.

Conversely, if we are given a monomial, then, by arranging its factors in a
suitable order, there is always a unique way to represent it as a two rowed
array satisfying the condition just given. The reader may check that one can
set up an \emph{insertion} algorithm exactly inverting the deletion procedure
above (what was deleted last, must be inserted first). In combinatorics one
most often uses standard bitableaux for the investigation of sequences (or two
row arrays). Then insertion is more important than deletion.

Since insertion and deletion are inverse operations, one obtains

\begin{theorem}
The map $\KRS$ is a bijection between the set of standard bitableaux on
$\{1,\dots,m\}\times\{1,\dots,n\}$ and the monomials of $K[X]$.
\end{theorem}

This theorem proves half of part (a) of the straightening law: it is enough to
check that every element of $K[X]$ can be written as a linear combination of
standard bitableaux. Using the straightening law one can now \emph{extend KRS
to a $K$-linear automorphism} of $K[X]$: with the standard representation
$x=\sum a_\Sigma\Sigma$ one sets $\KRS(x)=\sum a_\Sigma \KRS(\Sigma)$.

The automorphism KRS does not only preserve the total degree, but even the
$\ZZ^m\dirsum\ZZ^n$ degree introduced above: in fact, no column or row index
gets lost. Note also that KRS it is not a  $K$-algebra  isomorphism: it
acts as the identity on polynomials of degree $1$  but it is not the
identity map .  It would be interesting to have some insight on the
property of KRS as a linear map like, for instance, its eigenvalues and
eigenspaces.

\begin{remark}\label{KRSprop}
We note two important properties of KRS:
\begin{itemize}
\item[(a)] KRS commutes with
transposition of the matrix $X$: Let $X'$ be a $n\times m$ matrix of
indeterminates, and let
$\tau:K[X]\to K[X']$ denote the $K$-algebra  isomorphism induced by the
substitution
$X_{ij}\mapsto X'_{ji}$; then $\KRS(\tau(f))=\tau(\KRS(f))$ for all
$f \in K[X]$. Note that it suffices  to prove the equality   when $f$ is a
standard bitableau. Then the  statement  follows from \cite[Lemma 1.1]{HT}.

\item[(b)]
All the powers $\Sigma^k$ of a standard bitableau are again standard, and one
has
$$
\KRS(\Sigma^k)=\KRS(\Sigma)^k.
$$
\end{itemize}
\end{remark}

\section{\uppercase{KRS invariants and Gr\"obner bases}}\label{SectGr}

The power of KRS in the study of Gr\"obner bases for determinantal
ideals was detected by Sturmfels \cite{Stu}. He applied Schensted's theorem:

\begin{theorem}\label{Schensted}
Let $(t_1,\dots,t_w)$ be the shape of the standard bitableau $\Sigma$. Then
$t_1$ is the length of the longest strictly increasing subsequence in the lower
row of $\krs(\Sigma)$.
\end{theorem}

If $(v_{i_1},\dots,v_{i_q})$ is a strictly increasing subsequence of the lower
row, then the subsequence $(u_{i_1},\dots,u_{i_q})$ of the upper row must also
be strictly increasing. Therefore
\begin{equation}
\KRS(\Sigma)=M\cdot\diag[u_{i_1},\dots,u_{i_q}\sep
v_{i_1},\dots,v_{i_q}]\tag{$*$}
\end{equation}
where $M'$ is a monomial $\diag(\delta)$ denotes the product of all the
indeterminates in the diagonal of the  minor $\delta$.

Once and for all we now introduce a \emph{diagonal term order} on the
polynomial ring $K[X]$. With respect to such a term order the initial monomial
$\ini(\delta)$ is $\diag(\delta)$. There are various choices for a diagonal
term order, For example one can take the lexicographic order induced by the
total order of the $X_{ij}$ that coincides with the lexicographic order of the
$(i,j)$.

Schensted's theorem implies through its equivalent ($*$) that for a standard
bita\-bleau $\Sigma\in I_t$ there exists a $t$-minor $\delta$ such that
$$
\ini(\delta)=\diag(\delta) \mid \KRS(\Sigma),
$$
and, in particular, $\KRS(\Sigma)\in\ini(I_t)$: if $q>t$, then we can simply
write
$$
\diag[u_{i_1},\dots,u_{i_q}\sep v_{i_1},\dots,v_{i_q}] =
M''\diag[u_{i_1},\dots,u_{i_t}\sep v_{i_1},\dots,v_{i_t}].
$$
Since $I_t$ has a basis of standard bitableaux, it follows that
$\KRS(I_t)\subset \ini(I_t)$. The $K$-vector space $\KRS(I_t)$ has the same
Hilbert function as $I_t$ with respect to total degree since KRS preserves
total degree. But $\ini(I_t)$ also has the same Hilbert function as $I_t$. This
implies:

\begin{theorem}\label{Sturm}
The $t$-minors of $X$ form a Gr\"obner basis of $I_t$, and
$\KRS(I_t)=\ini(I_t)$.
\end{theorem}

In fact, the equation $\KRS(I_t)=\ini(I_t)$ has just been observed, and if
$M\in\ini(I_t)$ is a monomial, then it must be of the form $\KRS(\Sigma)$ for
some standard bitableau $\Sigma\in I_t$. But then $M$ is divisible by
$\ini(\delta)$ for some $t$-minor $\delta$. Exactly this condition must be
satisfied for the set of $t$-minors to form a Gr\"obner basis. It is worth
formulating the idea behind the proof of Theorem \ref{Sturm} as a lemma:

\begin{lemma}\label{KRS1}
\begin{itemize}
\item[(a)]
Let $I$ be an ideal of $K[X]$ which has a $K$-basis, say $B$, of standard
bitableaux, and let $S$ be a subset of $I$. Assume that for all $\Sigma\in B$
there exists $s\in S$ such that $\ini(s)\sep\KRS(\Sigma)$. Then $S$ is a
Gr\"obner basis of $I$ and $\ini(I)=\KRS(I)$.

\item[(b)]  Let $I$ and $J$ be  homogeneous ideals such that $\ini(I)=\KRS(I)$
and $\ini(J)=\KRS(J)$. Then $\ini(I)+\ini(J)=\ini(I+J)=\KRS(I+J)$ and
$\ini(I)\sect \ini(J)=\ini(I\sect J)=\KRS(I\sect J)$.
\end{itemize}
\end{lemma}

The proof of part (a) has been explained for the special case of $I=I_t$. For
(b) one uses
\begin{align*}
\KRS(I+J)&=\KRS(I)+\KRS(J)=\ini(I)+\ini(J)\subseteq \ini(I+J),\\
\KRS(I\sect J)&=\KRS(I)\sect \KRS(J)=\ini(I)\sect \ini(J)\supseteq
\ini(I\sect J),
\end{align*}
and concludes equality from the Hilbert function argument.

\begin{definition}
Let $I$ be an ideal with a standard basis. Then we say that $I$ is
\emph{in-KRS} if $\ini(I)=\KRS(I)$; if, in addition, the bitableaux $\Delta\in
I$ form a Gr\"obner basis, then $I$ is \emph{G-KRS}.
\end{definition}

In slightly different words, an ideal $I$ with a standard basis is in-KRS
if for each $\Sigma\in I$ there
exists $x\in I$ with $\KRS(\Sigma)=\ini(x)$; it is G-KRS if $x$ can always be
chosen as a bitableau.
As a consequence of Lemma \ref{KRS1} one obtains

\begin{lemma}\label{GandIn}
Let $I$ and $J$ be ideals with a basis of standard bitableaux.
\begin{itemize}
\item[(a)] If $I$ and $J$ are G-KRS, then $I+J$ is also G-KRS.
\item[(b)] If $I$ and $J$ are in-KRS, then $I+J$ and $I\cap J$ are also in-KRS.
\end{itemize}
\end{lemma}

In general the property of being G-KRS is not inherited by intersections as we
will see below.

The KRS correspondence could be used for the proof of Theorem \ref{Sturm}
since, by Schensted's theorem, the length of the first row of a standard
bitableau is a KRS invariant:

\begin{definition}
Let $\DD$ be the set of all bitableaux on the matrix $X$ and $F:\DD\to\NN$ a
function on $\DD$. Then we define a function on the set $\MM$ of monomials,
also called $F$, by
$$
F(M)=\max\{F(\Delta): \Delta\in\DD,\ M=\ini(\Delta)\}.
$$
Of course $F(M)$ is well-defined since there are only finitely many $
\Delta\in\DD$ with $M=\ini(\Delta)$.  We say $F$ is a \emph{KRS-invariant}
if the following conditions are satisfied:
\begin{itemize}
\item[(a)]
$F(\Delta)$ is the minimum of $F(\Sigma)$ where $\Sigma$ runs through the
standard bitableaux in the standard representation of a bitableau $\Delta$;
moreover, if $\Sigma'$ appears in the standard representation of $x\Delta$ for
some $x\in R$, then $F(\Sigma')\ge F(\Delta)$.
\item[(b)]
$
F(\Sigma)=F(\KRS(\Sigma))
$
for all standard bitableaux $\Sigma\in\DD$.
\end{itemize}
If just condition (a) is satisfied, then we say that we say that $F$ is
\emph{str-monotone}.
\end{definition}

In order to interpret condition (b) combinatorially we write
$$
\krs(\Sigma)=\begin{pmatrix}u_1&\dots&u_w\\v_1&\dots&v_w\end{pmatrix}.
$$
Then the bitableaux $\Delta$ such that $\KRS(\Sigma)=\ini(\Delta)$ correspond
bijectively to the decompositions of the lower row of $\krs(\Sigma)$ into
strictly increasing sequences, called \emph{inc-decompositions}. In fact, if
the sequence $v_{i_1},\dots,v_{i_t}$ is strictly increasing, then the same
holds for $u_{i_1},\dots,u_{i_t}$, or equivalently, $X_{u_1v_1}\cdots
X_{u_tv_t}$ is the diagonal product of a $t$-minor. Note that
$u_{i_1},\dots,u_{i_t}$ is always non-decreasing, and therefore
$X_{u_1v_1}\cdots X_{u_tv_t}$ can only be a diagonal product if
$v_{i_1},\dots,v_{i_t}$ is strictly increasing.

Thus condition (b) requires that $F(\Sigma)$, a number associated with the
standard bita\-bleau, is encoded in the sequence of integers forming the lower
row of $\krs(\Sigma)$.

It is now an easy exercise to show

\begin{proposition}\label{KRSinv}
Let $F$ be a KRS-invariant, and let  $k$ be an integer. Let $I_k(F)$ be the
ideal generated by all
bitableaux $\Delta$ such that $F(\Delta)\ge k$. Then
\begin{itemize}
\item[(a)] $I_k(F)$ has a standard  basis formed by all standard bitableaux
$\Sigma$ such that $F(\Sigma)\ge k$.
\item[(b)] Moreover,  $I_k(F)$ is G-KRS.
\end{itemize}
\end{proposition}

Part (a) follows immediately from str-monotonicity and   part (b) follows
from \ref{KRS1}. In general,  it does not suffice  to take the standard
bitableaux in  $I_k(F)$  to obtain a  Gr\"obner basis of $I_k(F)$;  we will
discuss an
example below.

Starting from the ideals $I_k(F)$ and applying  \ref{GandIn} one can now
find new ideals that are G-KRS or
at least in-KRS.

We have seen that the length of the first row is a KRS-invariant.
In order to apply KRS to a wider class of ideals one has to find other
(or more general) KRS invariants. One such family of invariants are the
functions
$\gamma_t$ defined as follows. For an integer $s$ and a sequence
$s_1,\dots,s_w$ of integers one sets
$$
\gamma_t(s)=(s-t+1)_+\qquad\text{and}\qquad \gamma_t(s_1,\dots,s_w)=
\sum_{i=1}^w \gamma_t(s_i).
$$
Here we have used the notation $(k)_+=\max(0,k)$. One then defines this
function for bitableaux $\Delta=\delta_1\cdots\delta_w$ by
$$
\gamma_t(\Delta)=\gamma_t(|\Delta|).
$$

The invariants $\gamma_t$ are of interest since they describe the symbolic
powers of the ideals $I_t$. Provided the characteristic of the field is $0$ or
$\ge \min(m,n)$ (we then say $K$ has \emph{non-exceptional} characteristic),
all products $I_{t_1}\cdots I_{t_r}$ have a primary decomposition as
intersections of such symbolic powers, and can therefore described in terms of
the $\gamma_t$:
\begin{itemize}
\item[(a)] $\Delta\in I_t^{(k)} \iff \gamma_t(\Delta)\ge k$;
\item[(b)] the standard bitableaux $\Sigma$ with $\gamma_t(\Sigma)\ge k$ are a
$K$-basis of $I_t^{(k)}$;
\item[(c)]
$$
I_{t_1}\cdots I_{t_r}= \Sect_{j=1}^s I_j^{(g_j)},\qquad
g_j=\gamma_j(t_1,\dots,t_r).
$$
\end{itemize}
See \cite[Section 10]{BV} and \cite{DEP}. The straightening law shows that
$\gamma_t$ is str-monotone. In \cite{BC} we have proved

\begin{theorem}\label{gammat}
The functions $\gamma_t$ are KRS-invariants.
\end{theorem}

As a consequence of this theorem and Lemma \ref{GandIn} one obtains that all
ideals $I_t^{(k)}$ are G-KRS and that all products of ideals of minors are
in-KRS if $\chara K=0$ or $\chara K\ge\min(m,n)$. Furthermore one can then show
that the ``initial algebras'' of the symbolic and ordinary Rees algebras of the
ideals $I_t$ are normal semigroup rings. In particular this implies that these
algebras are Cohen--Macaulay. Another object accessible to this approach is the
subalgebra $A_t$ of $K[X]$ generated by the $t$-minors of $X$. See \cite{BC}
for a detailed discussion.

\begin{examples}
(a) We choose $m,n\ge 3$. By the above discussion the ideal $I_2^{(2)}$ is
G-KRS. We want to show that the standard bitableaux in $I_2^{(2)}$ do not
form a Gr\"obner basis. The monomial $M=x_{12}x_{23}x_{21}x_{32}$ is the
initial term of a bitableau of shape $(2,2)$ and hence $M \in
\ini(I_2^{(2)})$.   If the standard bitableaux in $I_2^{(2)}$ were a
Gr\"obner basis, then $M$ would be divisible by the initial of a standard
bitableaux  in $I_2^{(2)}$. The   standard bitableaux in $I_2^{(2)}$ of
degree $\leq 4$ have shape  $(3)$, $(3,1)$ and $(2,2)$ and clearly their
initial term cannot divide $M$.

(b) Suppose that $I$ and $J$ are G-KRS and let $\Sigma$ be a standard
bitableau
in $I\cap J$. Then we can find bitableaux $\Delta_1\in I$ and $\Delta_2\in J$
such that $\KRS(\Sigma)=\ini(\Delta_1)=\ini(\Delta_2)$. In general it can
happen that $\Delta_1\notin J$ and $\Delta_2\notin I$, and $I\cap J$ need not
be G-KRS. An example is $I=I_3^{(2)}$, $J=I_4$ for a matrix of size at least
$6\times 6$. In fact, let $\Sigma=[1\,3\,4\,5\sep 1\,2\,3\,6][2\,6\sep 4\,5]$.
(This is the example considered in Section \ref{SectKRS}). Then $\Sigma\in J$
is obvious, and $\Sigma\in I$ since $\gamma_2(\Sigma)\ge 2$. Since $I\cap J$ is
in-KRS, it follows that $\KRS(\Sigma)=
X_{14}X_{21}X_{32}X_{45}X_{56}X_{}X_{63}\in \ini(I\cap J)$. It is however
impossible to write this monomial as the initial monomial of a bitableau
$\Delta$ in $I\cap J$. The bitableaux of  degree $6$  in $I\cap J$ are those of
shapes $(6)$, $(5,1)$, and $(4,2)$. (By Schensted's theorem, only the last
shape would be possible.)
\end{examples}

Let $T$ be a new indeterminate and define the ideal $J$ in the extended
polynomial ring $K[X][T]$ by
$$
J=I_m+I_{m-1}T+\dots+I_1T^{m-1}+(T^m);
$$
we assume that $m\le n$. This ideal and its Rees algebra $\Dirsum_{i=0}^\infty
J^i$ is fundamental for the generic case of MacPherson's graph construction;
see \cite{BK}. If one ``expands'' the power $J^k$ into a ``polynomial'' in $T$,
then the ``coefficient'' of $T^{km-d}$ is
\begin{multline*}
J(k,d)=\sum I_0^{e_0}I_1^{e_1}\cdots I_m^{e_m},\\ e_0+e_1+\dots+e_m=k,\
e_1+2e_2+\dots+me_m=d.
\end{multline*}
For a non-increasing sequence $s_1,\dots,s_w$ of non-negative integers let us
define
$$
\alpha_k(s_1,\dots,s_w)=\sum_{i=1}^k s_i
$$
where $s_i=0$ if $i>w$. Then we can set
$$
\alpha_k(\Delta)=\alpha_k(|\Delta|)
$$
for every bitableau $\Delta$. The straightening law shows that $\alpha_k$ is
str-monotone, and it follows easily that
\begin{itemize}
\item[(a)] $\Delta\in J(k,d) \iff \alpha_k(\Delta)\ge d$;
\item[(b)] the standard bitableaux $\Sigma$ with $\alpha_k(\Sigma)\ge d$ are a
$K$-basis of $J(k,d)$.
\end{itemize}

\begin{theorem}\label{Greene}
The functions $\alpha_k$ are KRS-invariants.
\end{theorem}

An analysis of $\alpha_k$ in terms of inc-decompositions shows that
$\alpha_k(KRS(\Sigma))\ge d$ if and only if the lower sequence of
$\krs(\Sigma)$ contains a subsequence of length $d$ that itself can be
decomposed into $k$ increasing subsequences. Thus Theorem \ref{Greene} is just
a re-interpretation of Greene's theorem \cite{Gr}: $\alpha_k(\Sigma)$  is the
maximal length of a subsequence that has an inc-decomposition into $k$ parts.

For the ``determinantal'' consequences of Theorem \ref{Greene} we refer the
reader to \cite{BK}. The relationship between the KRS invariants $\gamma_t$ and
$\alpha_k$ is analyzed in Section \ref{SectBas}.

\section{\uppercase{Ideals defined by shape}}\label{SectSh}

We say that an ideal $I\subset K[X]$ is \emph{defined by shape} if it is
generated as an ideal by a set of bitableaux, and, moreover, it depends
only on $|\Delta|$
whether a bitableau $\Delta$ belongs to $I$. In this section we want to
characterize the ideals defined by shape in the case in which the
characteristic of $K$ is big enough. In particular we will see that all these
ideals are in-KRS.

The following ``balancing lemma'' is a crucial argument; it is a simplified
version of \cite[(10.10)]{BV}.

\begin{lemma}\label{bal}
Let $\pi$ and $\rho$ be minors of $X$, and set $u=|\rho|$, $v=|\pi|$ (we
include the case $\pi=1$, in which $u=0$). Suppose that $u<v$ and $\chara K=0$
or $\chara K>\min\bigl(u+1,m-(u+1),n-(u+1)\bigr)$. Then $\pi\rho\in
I_{u+1}I_{v-1}$.
\end{lemma}

The case in which $u=0$ is just Laplace expansion. In general the lemma says
that a product $\pi\rho$ of minors can be expressed as a linear combination of
minors that are ``more balanced'' in size. By repeated application of the lemma
we see that $\pi\rho$ is even a linear combination of products $\delta\epsilon$
such that $|\delta|+|\epsilon|=|\pi|+|\rho|$ and $|\delta| \le |\epsilon|\le
|\delta|+1$.

The group $\GL=\GL(m,K)\times \GL(n,K)$ operates as a group of linear
substitutions on $R=K[X]$ in a natural way: for $M\in\GL(m,K)$ and
$N\in\GL(n,K)$ one substitutes $X_{ij}$ by the corresponding entry of
$MXN^{-1}$. Therefore $R$ is an interesting object for the representation
theory of $\GL$, and representation theory offers another approach to the
theory of determinantal rings. It is clear that an ideal defined by shape is
$\GL$-stable since each element $g\in\GL$ transforms a minor into a linear
combination of minors of the same size.

Let $\sigma$ be a shape. Then the ideal $I^{(\sigma)}$ is generated by all
bitableaux $\Delta$ with $\gamma_t(\Delta)\ge \gamma_t(\sigma)$ for all $t$
(see \cite[Section 11]{BV}). If all these inequalities hold, then
$$
|\Delta|\ge\sigma.
$$
If even $|\Delta|_k\ge \sigma_k$ for all $k$, then we write $|\Delta|\supset
\sigma$.

The ideals $I^{(\sigma)}$ are evidently defined by shape. By definition
$$
I^{(\sigma)}=\Sect_{t} I_t^{(\gamma_t(\sigma))}.
$$
Therefore $I^{(\sigma)}$ has a standard  basis and is in-KRS.

The following theorem should have been contained in \cite[Section 11]{BV}.

\begin{theorem}\label{shape}
Suppose that $\chara K=0$. Then the following are equivalent for an ideal
$I\subset K[X]$:
\begin{itemize}
\item[(a)] $I$ is defined by shape;
\item[(b)] $I$ is a sum of ideals of type $I^{(\sigma)}$.
\item[(c)] $I$ has a standard basis, and $I$ is stable under the action of
$\GL$ on $K[X]$.
\end{itemize}
\end{theorem}

\begin{proof}
(a)$\implies$(b): $I$ is obviously contained in the sum of all ideals
$I^{(|\Delta|)}$ where $\Delta$ runs through the generators of $I$.

On the other hand, let $\Sigma$ be a (standard) bitableau contained in
$I^{(|\Delta|)}$, that is, $|\Sigma|\ge|\Delta|)$. If even
$|\Sigma|\supset|\Delta|$, then we can apply Laplace expansion and write
$\Sigma$ as a linear combination of bitableaux of the the same shape as
$\Delta$. Suppose that $|\Sigma|\not\supset|\Delta|$, and let $k$ be the
smallest index such that $|\Sigma|_k<|\Delta|_k$. Then there must be an index
$j<k$ such that $|\Sigma|_j>|\Delta|_j$. Now one applies the balancing lemma
above, increasing the $k$-th row at the expense of the $j$-th. After finitely
many balancing steps we have written $\Sigma$ as a $K$-linear combination of
bitableaux $\Xi$ such that $|\Xi|\supset|\Delta|$.

(b)$\implies$(a): This is evident, as well as (b)$\implies$(c).

(c)$\implies$(b): We have to show that $\Xi\in I$ for all standard bitableaux
$\Xi$ such that $|\Xi|\ge|\Sigma|$ for some standard bitableau $\Sigma\in I$.

Set $\sigma=|\Sigma|$. Suppose first that $|\Xi|=|\Sigma|$. By
\cite[(11.10)]{BV} there is a decomposition of $\GL$ stable subspaces
$$
I^{(\sigma)}=M_\sigma \dirsum I_>^{(\sigma)}
$$
where $I_>^{(\sigma)}$ is generated by all standard bitableaux $\Theta$ such
that $|\Theta| > \sigma$; in this decomposition $M_\sigma$ is irreducible and
the unique $\GL$-stable complement of $I_>^{(\sigma)}$ in $I^{(\sigma)}$. Thus
we can write
$$
\Xi=x+y
$$
where $x\in M_\sigma$, $y\in I_>^{(\sigma)}$. Since $\Xi\notin I_>^{(\sigma)}$,
we deduce that $x\neq 0$.

For $\Xi=\Sigma$ it follows that $M_\sigma\subset I$. In fact, $I$ has a unique
decomposition as a direct sum of irreducible $\GL$-modules, and these are
exactly the $M_\tau$; since the projection to $M_\sigma$ is non-trivial, it
must appear in the decomposition of $I$.

For general $\Xi$ it now follows that $x\in I$, and therefore all the standard
bitableaux in the standard representation of $x$ must belong to $I$. However,
since $y\in I_>^{(\sigma)}$, $\Xi$ must appear in this standard representation.

Suppose now that $|\Xi|>\sigma$. If we find some standard bitableau
$\Gamma$ of shape $|\Xi|$ such that
$\Gamma\in I$, then the argument just given shows that $\Xi\in I$. Since
we can connect  $\sigma$ and $|\Xi|$ by a
chain of shapes with respect to the partial order $\le$, it is enough to
consider the case in which $|\Xi|$ is an upper neighbor of $\sigma$. One
obtains the upper neighbors by either inserting a new box below the bottom of
the diagram (and thereby increasing $\gamma_1$ by $1$) or by removing an
``outer corner'' box of the Young diagram of shape $\sigma$ and inserting it at
an ``inner corner'' in such a way that the box travels one row up or one column
to the right. (At the end of the first row is also an inner corner, provided
its length is $<\min(m,n)$). Figure \ref{Upper} illustrates the three cases.
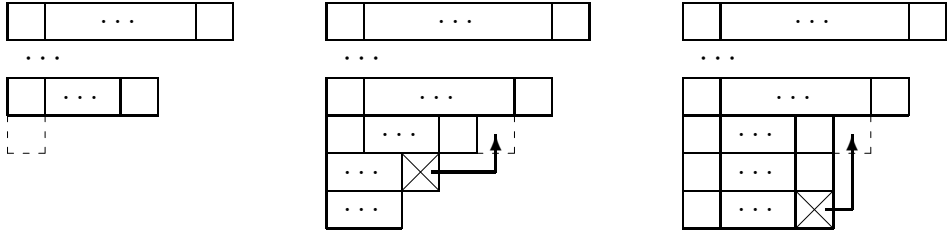
\begin{figure}[hbt]
$$
\unitlength=0.5cm
\begin{picture}(6,6)(0,0)
\EBox05\LBox15{\dots}4\EBox55
\EBox03\LBox13{\dots}2\EBox33
\EDBox02
\put(0,4){\makebox(2,1){$\cdots$}}
\end{picture}
\qquad\quad
\begin{picture}(7,6)(0,0)
\EBox05\LBox15{\dots}5\EBox65
\EBox03\LBox13{\dots}4\EBox53
\EBox02\LBox12{\dots}2\EBox32\EDBox42
\LBox01{\dots}2\EBox21\Cross21
\LBox00{\dots}2
\put(0,4){\makebox(2,1){$\cdots$}}
\thicklines
\put(2.8,1.5){\line(1,0){1.7}}
\put(4.5,1.5){\vector(0,1){1.0}}
\end{picture}
\qquad\quad
\begin{picture}(7,6)(0,0)
\EBox05\LBox15{\dots}5\EBox65
\EBox03\LBox13{\dots}4\EBox53
\EBox02\LBox12{\dots}2\EBox32\EDBox42
\EBox01\LBox11{\dots}2\EBox31
\EBox00\LBox10{\dots}2\EBox30\Cross30
\put(0,4){\makebox(2,1){$\cdots$}}
\thicklines
\put(3.8,0.5){\line(1,0){0.7}}
\put(4.5,0.5){\vector(0,1){2.0}}
\end{picture}
$$
\caption{Upper neighbors with respect to $\le$}\label{Upper}
\end{figure}

The case in which a new box is added, is trivial. In fact, we fill it with
$X_{mn}$, and $X_{mn}\Sigma$ is standard of the right shape.

We now assume that a box travels one row up, say from row $k$ to row $k-1$.
With $\sigma=(s_1,\dots,s_w)$ let $p=s_k$ and $q=s_{k-1}$. One forms the
standard bitableau $\Theta=\theta_1\cdots\theta_w$ where
\begin{align*}
\theta_j&=[1,\dots, s_j\sep 1,\dots, s_j],\qquad j\neq k,\\ \theta_k&=[1,\dots,
p-1, q+1\sep 1,\dots, p-1, q+1]
\end{align*}
of shape $\sigma$. Since $\Theta$ is a standard bitableau of
shape $\sigma$,  it belongs to $I$, as was shown above. Note that
$s_{k+1}<p<q<s_{k-2}$ (where $s_{k+1}=0$ if $k=w$ and $s_{k-2}=\infty$ if
$k=2$).

Now we apply the cyclic permutation $\pi\in (p\ p+1\dots q\ q+1)$ to both the
rows and columns of the matrix $X$; the transformation $\pi$ belongs to $\GL$.
Therefore $\pi(\Theta)\in I$. All the factors of $\Theta$ except $\theta_{k-1}$
and $\theta_k$ are invariant under $\pi$, whereas
\begin{align}
\pi(\theta_{k-1})&=[1,\dots,p-1,p+1,\dots,q+1\sep 1,\dots,p-1,p+1,\dots,q+1]\\
\pi(\theta_{k})&=[1,\dots,p-1,p\sep 1,\dots,p-1,p].
\end{align}
If we straighten this product and multiply its standard presentation with the
remaining factors of $\Theta$, then we obtain the standard representation of
$\pi(\Theta)$. Therefore it is enough that a standard bitableau of shape
$(q+1,p-1)$ appears in the standard representation of
$\pi(\theta_{k-1})\pi(\theta_k)$. (Whereas the proof of ((a)$\implies$(b) is
based on ``balancing'', we now need the ``unbalancing'' effect of
straightening.) Mapping all the indeterminates $X_{ij}$ with $i\neq j$, $i<p$
or $j<p$, to $0$ and $X_{11},\dots,X_{p-1\, p-1}$ to $1$, one reduces the claim
to the assertion that in the standard representation of $[1\sep
1][2,\dots,r\sep 2,\dots, r]$ with $r\ge 3$ the minor $[1,\dots,r\sep
1,\dots,r]$ shows up, and that is immediate from Laplace expansion.

The case in which the box travels one column to the right is similar and left
to the reader. Essentially it is the case in which $\sigma$ consists of a
single column.
\end{proof}

One should note that the implications (b)$\implies$(a) and (b)$\implies$(c) are
true over arbitrary fields (actually, over all rings of coefficients), whereas
(a)$\implies$(b) needs only that $\chara K>\min(m,n)$. The implication
(c)$\implies$(b) uses the hypothesis that $K$ has characteristic $0$ more
profoundly: the ideal $I_1^{p+1}+(X_{11}^p,\dots,X_{mn}^p)$ satisfies (c) if
$\chara K=p>0$, but is not a sum of ideals $I^{(\sigma)}$ (provided that $X$ is
not just a $1\times1$ matrix).

\begin{corollary}\label{AllIn}
Suppose that $K$ is a field of characteristic $0$. Then all $\GL$-stable ideals
that have a standard basis are in-KRS.
\end{corollary}

This follows from Lemma \ref{GandIn} since the ideals $I^{(\sigma)}$ are
in-KRS.

Note that there are ideals with a standard basis which are not in-KRS.  For
instance, let  $\Sigma$  be   standard bitableau  with
$\KRS(\Sigma)\neq \ini(\Sigma)$ and set $d=\deg\Sigma$. Then
$I=(\Sigma)+(X_{ij} : 1\leq i\leq m, 1\leq j \leq n)^{d+1}$  has a standard
basis and $\KRS(I)\neq \ini(I)$.
As we have seen in the previous section, there are in-KRS ideals that are not
G-KRS. The G-KRS ideals among those considered in Corollary \ref{AllIn} will be
characterized in the next section.

\section{\uppercase{Basic KRS-invariants}}\label{SectBas}

In this section we want to show that the functions $\alpha_k$ are basic
KRS-invariants, as far as functions $F:\DD\to\NN$ are considered that depend
only on shape. First we prove a converse of Proposition \ref{KRSinv}:

\begin{proposition}\label{InvConv}
Let $F:\DD\to\NN$ be a str-monotone function, and let $k$ be an integer.
Let  $I_k(F)$  be the ideal generated by all bitableaux $\Delta$ such that
$F(\Delta)\ge k$. Then
the following hold:
\begin{itemize}
\item[(a)] $I_k(F)$ has a standard basis, and $\Delta\in I_k(F)$ (if and)
only if $F(\Delta)\ge k$;
\item[(b)] if $I_k(F)$ is G-KRS for all $k$, then $F$ is a KRS-invariant.
\end{itemize}
\end{proposition}

\begin{proof}
(a) It follows immediately from the definition of str-monotonicity that
$I_k(F)$ is the $K$-vector space generated by the bitableaux $\Delta$ with
$F(\Delta)\ge k$, and that the standard bitableaux $\Sigma\in I_k(F)$ form a
$K$-basis.

(b) We have to show that $F(\Sigma)=F(\KRS(\Sigma))$ for every standard
bitableau $\Sigma$. First note that  $F(\Sigma)\leq F(\KRS(\Sigma))$. Set
$k=\KRS(\Sigma)$. Since by assumption $I_k(F)$ is G-KRS we have that
$\KRS(\Sigma)=\ini(\Delta)$ for some $\Delta\in I_k(F)$ and consequently
$F(\KRS(\Sigma))\ge k$. It follows that for every $t$  one has
$$
\KRS(I_t(F)) \subseteq \langle M : F(M)\geq t \rangle  \subseteq
\ini(I_t(F)) =\KRS(I_t(F))
$$
and hence
$$\KRS(I_t(F)) = \langle M : F(M)\geq t \rangle.$$
Now let $M$ be a monomial and set $t=F(M)$. Then
$\Sigma=\KRS^{-1}(M)$ is  in $I_t(F)$ which proves that
$F(\Sigma)\geq F(\KRS(\Sigma))$.
\end{proof}

Now we can show that the KRS-invariance of the functions $\gamma_t$ follows
from that of the $\alpha_k$:

\begin{proposition}\label{algam}
For all $t$, $1\le t\le \min(m,n)$ and all $r\ge1$ one has
$$
I_t^{(r)}=\sum_{k\ge 1} J\bigl(k,r+k(t-1)\bigr).
$$
\end{proposition}
\begin{proof}
Let $\Sigma\in I_t^{(r)}$ be a (standard) bitableau, and suppose that $k$ is
the biggest index such that $|\Sigma|_k\ge t$. Then obviously $\Sigma\in
J(k,r+(t-1))$.

The verification of the inclusion $\supset$ is likewise simple.
\end{proof}

Since by Greene's theorem the ideals $J(k,d)$ are G-KRS, it follows immediately
that the ideals $I_t^{(r)}$ are G-KRS. In conjunction with Proposition
\ref{InvConv} we therefore obtain that the $\gamma_t$ are KRS-invariants. From
hindsight, Theorem \ref{gammat} is an easy consequence of Greene's theorem.

The most difficult part of \cite{BC} is the proof that the ideals $I_t^r$ are
G-KRS in non-excep\-tion\-al characteristics. Actually $I_t^r=J(r,rt)$ so that the
property of being G-KRS is no longer surprising for $I_t^r$.

Not every shape defined, G-KRS ideal is the sum of ideals $J(k,d)$. As an
example one can take $I_1^4\cap I_3=I_1I_3$. Another exception occurs for
square matrices. For example, if $m=n=4$, then $I_2I_4$ is G-KRS, but it looses
this property for bigger matrices. Nevertheless, the next theorem shows that
the functions $\alpha_k$ are truly basic KRS-invariants:

\begin{theorem}
Let $I$ be a shape defined ideal. If $I$ is G-KRS, then it is the sum of ideals
of type $J(k,d)\cap I_1^u$ and, if $m=n$, $\bigl(J(k,d)\cap I_1^u\bigr)I_n^v$.
\end{theorem}

\begin{proof}
Let $\Sigma\in I$ be a (standard) bitableau. Suppose that
$$
\sigma=|\Sigma|=(s_1,\dots,s_k,1,\dots,1)
$$
with $s_k\ge 2$, and $s_1<\max(m,n)$ (equality can only occur if $m=n$). Set
$d=\alpha_k(\sigma)$ and $D=\deg(\Sigma)$. Then $J=I(k,d)\cap I_1^D$ is the
smallest ideal of type $J(k,d)\cap I_1^u$ containing $\Sigma$.

Among the shapes of the elements in the standard basis of $J$ there exists a
unique element that is minimal with respect to the partial order defined by the
functions $\gamma_t$ (see Section \ref{SectSh}). In fact let $u=\lfloor
d/k\rfloor$ and $r=d-uk$. Then the smallest element is
$$
\theta=(u+1,\dots,u+1,u,\dots,u,1\dots,1)
$$
where $u+1$ appears $r$ times, $u$ appears $d-r$ times, and $1$ as often as in
$\sigma$. We have $J=I^{(\theta)}$. Therefore, if $\sigma=\theta$, the ideal
$J$ is contained in $I$, as follows from Theorem \ref{shape}.

Now suppose that $\sigma>\theta$. Then it is enough to show that there exists a
standard bitableau $\Sigma$ of shape $\sigma$ such that $|\Delta|<\sigma$
whenever $\ini(\Delta)=\KRS(\Sigma)$. Since we have $|\Delta|\le \sigma$ from
the KRS invariance of the functions $\gamma_t$ (or $\alpha_k$) it really
suffices to find $\Sigma$ of shape $\sigma$ such that $|\Delta|\neq \sigma$ for
all $\Delta$ with $\ini(\delta)=\KRS(\Sigma)$.

Instead of constructing $\Delta$ directly, we find a monomial $M$ such that
$\KRS^{-1}(M)$ has the desired shape, but $M$ cannot be written as
$\ini(\Delta)$ where $|\Delta|=\sigma$. The shape of $\KRS^{-1}(M)$ can be
controlled via the $\gamma$ or $\alpha$ functions.

Let us first consider $\sigma=(s_1,s_2)$; we set $p=s_1$, $q=s_2$. Then
$\max(m,n)>p\ge q+2>1$ and $\min(m,n)\ge p$. It is harmless to assume
$n=\max(m,n)$. With $r=p-q+2$ let
$$
M=X_{11}X_{22}X_{34}\cdots X_{p,p+1}\cdot X_{13}\cdot X_{r3}X_{r+1,4}\cdots
X_{p,q+1}.
$$
For $p=5$ and $q=3$ this monomial has the following ``picture'':
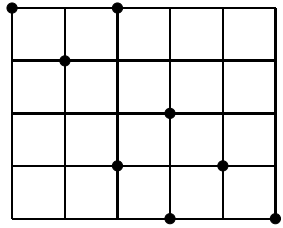
$$
\begin{picture}(5,4)(0,0)
 \def\vertex{\circle*{0.20}}
 \multiput(0,0)(1,0){6}{\line(0,1)4}
 \multiput(0,0)(0,1){5}{\line(1,0)5}
 \put(0,4){\vertex}
 \put(1,3){\vertex}
 \put(3,2){\vertex}
 \put(4,1){\vertex}
 \put(5,0){\vertex}
 \put(2,4){\vertex}
 \put(2,1){\vertex}
 \put(3,0){\vertex}
\end{picture}
$$

The reader can check that $\KRS^{-1}(M)$ indeed as shape $(p,q)$. However, it
is not possible to decompose $M$ into a product $\ini(\delta_1)\ini(\delta_2)$
where $|\delta_1|=p$, $|\delta_2|=q$.

In the general case one must multiply $M$ with suitable factors. These are not
hard to find.
\end{proof}

\section{\uppercase{Initial ideal of symbolic powers and symbolic powers of the
initial ideal}}\label{SectIn}

Let $S$ be a polynomial ring and $\tau$ a term order on $S$. Given an ideal $J$
of $S$ and an integer $b$  one denotes by $J_{\leq b}$ the intersection of all
the primary components of height $\leq b$ in a primary decomposition of $J$
($J_{\leq b}$ is independent of the chosen primary decomposition). One knows
that $J_{\leq b}=\{ f : \height (J:f)>b\}$ and that
$$
\ini(J_{\leq b})\subseteq \ini(J)_{\leq b}, \eqno{(1)}
$$
see \cite{STV}. If $J$ has height $c$, then we define the $k$-symbolic power
$J^{(k)}$ of $J$ to be:
$$
J^{(k)}=(J^k)_{\leq c}.
$$
Note that $J^{(1)}=J$ if and only if $J$ has no embedded primes and all the
minimal primes have height $c$. Now by $(1)$
$$
\ini(J^{(k)})\subseteq \ini(J^k)_{\leq c}.
$$
On the other hand, since $\ini(J)^k\subseteq  \ini(J^k)$,  we have
$$
\ini(J)^{(k)}=(\ini(J)^k)_{\leq c}  \subseteq \ini(J^k)_{\leq c}.
$$
Summing up, there are inclusions
$$
\ini(J^{(k)})\subseteq \ini(J^k)_{\leq c} \supseteq \ini(J)^{(k)}. \eqno{(2)}
$$
These inclusions are in general strict. In this section we will show that they
are equalities when $J$ is the determinantal ideal $I_t$ and $\tau$ is a
diagonal term order.

\begin{theorem}\label{symb-in}
In non-exceptional characteristics we have
$$
\ini(I_t^{(k)})= \ini(I_t^k)_{\leq c} = \ini(I_t)^{(k)}
$$
where $c$ is the height of $I_t$, i.e. $c=(m-t+1)(n-t+1)$.
\end{theorem}

The initial ideal $\ini(I_t)$ of $I_t$ is the square free monomial ideal
generated by the diagonals of the $t$-minors, simply called $t$-diagonals in
the sequel. Hence it is the Stanley-Reisner ideal of the simplicial complex
$$
\Delta_t=\{ A\subseteq \{1,\dots,m\}\times \{1,\dots,n\}: A \mbox{ does
not contain $t$-diagonals}\}.
$$
Denote by  ${\bf F}_t$  the set of the facets of $\Delta_t$. Then
$$
\ini(I_t)=\Sect_{F \in {\bf F}_t} P_F
$$
where $P_F$ denotes the ideal generated by the $x_ {ij}$ with $(i,j) \not\in
F$. The elements of ${\bf F}_t$  are described in \cite{HT} in terms of
families of non-intersecting paths. It turns out that $\Delta_t$ is a pure
(even shellable) simplicial complex. Since the powers of the ideals $P_F$ are
$P_F$-primary, it follows that
$$
\ini(I_t)^{(k)}=\Sect_{F \in {\bf F}_t} P_F^k.
$$
We start by proving:

\begin{proposition}\label{symb-in1}
We have $$\ini(I_t^{(k)})=\Sect_{F \in {\bf F}_t} P_F^k.$$
\end{proposition}

We have already mentioned in Section \ref{SectGr} that  $I_t^{(k)}$ is G-KRS,
and in particular the initial ideal $\ini(I_t^{(k)})$ of $I_t^{(k)}$ is
generated by the monomials $M$ with $\gamma_t(M)\geq k$. Now a monomial
$M=\prod_{i=1}^s x_{a_ib_i}$ is in $P_F^{k}$ if and only if the cardinality of
$\{i : (a_i,b_i) \not\in F\}$ is $\geq k$. Equivalently,  $M$ is in $P_F^{k}$
if and only if the cardinality of $\{i : (a_i,b_i) \in F\}$ is $\leq \deg(M) -
k$.  If we set
$$
\w_t(M)=\max\bigl\{|A| :  A\subseteq [1,\dots,s] \mbox{ and }  \{(a_i,b_i) : i
\in A\}\in \Delta_t\bigr\}
$$
then we have that a monomial $M$ is in $\Sect_{F\in {\bf F}_t} P_F^{k}$ if and
only if $\w_t(M)\leq \deg(M) -k $, that is, $ \deg(M)-\w_t(M) \geq k$. Now
Proposition \ref{symb-in1} follows from:

\begin{lemma}\label{symb-in2}
Let $M$ be a monomial.   Then $\gamma_t(M)+\w_t(M)=\deg(M)$.
\end{lemma}

We reduce this lemma to a combinatorial statement on sequences of integers.
Given such a sequence $b$ we define $\w_t(b)$ to be the cardinality of the
longest subsequence of $b$ which does not contain an increasing subsequence of
length $t$, that is,
$$
w_t(b)=\max\{ \length(c) : c \mbox{ is a subsequence of $B$ and }
\gamma_t(c)=0\}.
$$
Let  $M=\prod_{i=1}^s x_{a_ib_i}$ be a monomial. We may order the indices such
that $a_i\leq a_{i+1}$ for every $i$ and $b_{i+1}\geq b_{i}$ whenever
$a_i=a_{i+1}$. (We have already considered this rearrangement in Section
\ref{SectKRS}.) Then the $t$-diagonals dividing $M$ correspond to increasing
subsequences of length $t$ of the sequence $b$, and $\w_t(M)=\w_t(b)$.  Since
$\w_t(M)$ depends only on the sequence $b$ we may assume that $a_i=i$ for every
$i$. Then, by exchanging the role between the $a_i$'s and the $b_i$'s we may
also assume that the $b_i$ are distinct integers (see Remark \ref{KRSprop}).
Summing up, it suffices to show that:

\begin{lemma}\label{symb-in3}
One has $\gamma_t(b)+\w_t(b)=\length(b)$ for every sequence $b$ of distinct
integers.
\end{lemma}

\begin{proof}
Let  $P$ be the tableau obtained from $b$ by the Robinson-Schensted insertion
algorithm. We have already discussed the first part of Greene's theorem, namely
that the sum $\alpha_k(P)$ of the lengths of the first $k$ rows of $P$ is the
length of the longest subsequence of $b$ that has a decomposition into $k$
increasing subsequences. But the theorem contains a second (dual) assertion:
the sum $\alpha_k^*(P)$ of the lengths of the first $k$ columns of $P$ is the
length of the longest subsequence of $b$ that can be decomposed into $k$
decreasing subsequences.

It follows that a sequence $a$ has no increasing subsequence of length $t$ if
and only if it can be decomposed into $t-1$ decreasing subsequences. Then
$w_t(b)$ is the equal to the maximal length of a subsequence of $b$ which can
be decomposed into $t-1$ decreasing subsequences.

Therefore $\w_t(b)=\alpha_{t-1}^*(P)$. On the other hand, by Theorem
\ref{gammat} we know that $\gamma_t(b)$ is equal to $\gamma_t(P)$ which is the
sum of the length of the columns of $P$ of index $\geq t$. Therefore
$\gamma_t(b)+\w_t(b)$ is equal to the number of entries of $P$ which is the
length of $b$.
\end{proof}

We know \cite[Thm.\ 3.5]{BC} that $\ini(I_t^k)=\Sect_{j=1}^t
\ini(I_j^{(k(t+1-j))})$ (in non-exceptional characteristic) and hence, taking
into consideration Proposition \ref{symb-in1},  we have:
$$
\ini(I_t^k)=\Sect_{j=1}^t \Sect_{F\in  {\bf F}_j} P_F^{k(t+1-j)}.
\eqno{(3)}
$$
Since the powers of the ideal $P_F$ are $P_F$-primary we have that
$(3)$ is indeed a primary decomposition of $\ini(I_t^k)$. Hence
$\ini(I_t^k)_{\leq c}$ is equal to $\Sect_{F\in  {\bf F}_t} P_F^{k}$.  This
concludes the proof of \ref{symb-in}.

\begin{remark}
(a) Theorem \ref{symb-in} can be interpreted as a description of the normal
semigroup of monomials in the initial algebra of the Rees algebra of $I_t$
in terms of linear inequalities.

(b)  The argument above shows also that $\ini((I_t^k)_{\leq
b})=\ini(I_t^k)_{\leq b}$ for every integer $b$. But in general
$(\ini(I_t)^k)_{\leq b}$ is strictly smaller than $\ini(I_t^k)_{\leq b}$. This
is because  not all the variables appear in the generators of $\ini(I_t)$ while
the maximal ideal is associated to $I_t^k$ (and hence to $\ini(I_t^k)$) for
large  $k$.
\end{remark}

\begin{question}
What is a primary decomposition of the powers of $\ini(I_t)$? Is the Rees
algebra of $\ini(I_t)$ Cohen-Macaulay? Is it normal?
\end{question}

\section{\uppercase{Cogenerated ideals}}\label{SectCo}

As before, let $X=(x_{ij})$ be an $m\times n$ matrix of indeterminates. We
consider the set of minors of $X$ equipped with the usual partial order that
has been introduced in Section \ref{SectStr}. Let
$\delta=[a_1,a_2,\dots,a_r|b_1,b_2,\dots,b_r]$ be a minor of $X$. One defines
$I(\delta,X)$ to be the ideal of $K[X]$ generated by all the minors $\mu$ such
that $\mu\not \succeq \delta$, i.e.
$$
I(\delta,X)=(\mu : \mu\not \succeq \delta).
$$

The ideal $I(\delta,X)$ is said to be the ideal {\emph cogenerated}  by
$\delta$. For general facts about  the ideals cogenerated by minors we refer
the reader to \cite{BV}. We just recall that $I(\delta,X)$ is a prime ideal and
that it has a standard basis. Namely the set
$$
B(\delta)=\{ \Sigma :  \Sigma=\sigma_1\cdots\sigma_w \text{ is a standard
bitableau and } \sigma_1\not \succeq \delta\}
$$
is  a $K$-vector space basis of $I(\delta,X)$. Herzog and Trung  have shown in
\cite[Theorem 2.4]{HT} that the natural generators of $I(\delta,X)$
(i.e.~the minors
$\mu$ such that $\mu\not \succeq \delta$) form a Gr\"obner basis of
$I(\delta,X)$ with respect to the diagonal term order. Their argument makes use
of the KRS correspondence and boils down to the study of the  KRS image of the
elements of $B(\delta)$.  We will see  that this can be rephrased in terms of
properties of a suitable $\gamma$-function associated to $\delta$. We will
henceforth denote its value on the bitableau $\Delta$ by
$\gamma_\delta(\Delta)$. We start by defining $\gamma_\delta(\mu)$ for a single
minor $\mu=[c_1,\dots,c_k\sep d_1,\dots,d_k]$, namely
$$
\gamma_\delta( \mu)= \max\bigl\{(i-j+1)_+ : 1\leq i\leq k, 1\leq j\leq r+1
\mbox{ and } (c_i< a_j \mbox{ or } d_i < b_j) \bigr\}
$$
where, by definition, $a_{r+1}=b_{r+1}=\infty$. Then we extend, by linearity,
the $\gamma_\delta$-function  to product of minors, that is, if
$\Delta=\mu_1\cdots \mu_h$ is a product of minors, then
$$
\gamma_\delta( \Delta)=\sum_{i=1}^h \gamma_\delta( \mu_i).
$$
The standard bitableaux $\Sigma$ such that $\gamma_\delta( \Sigma)\neq 0$ are
exactly the elements of $B(\delta)$. Note that if one takes
$\delta=[1,\dots,t-1\sep 1,\dots,t-1]$ then $I(\delta,X)=I_t$ and
$\gamma_\delta(\mu)=\gamma_t(\mu)$. We may extend, as we have done in Section
\ref{SectGr}, the definition of the $\gamma_\delta$-function also to ordinary
monomials by setting:
$$
\gamma_\delta(M)=\max\{\gamma_\delta(\Delta) : \Delta \mbox{ is a bitableau and
} \ini(\Delta)=M\}.
$$
In terms of the $\gamma_\delta$-function Herzog and Trung \cite[Lemma 1.2]{HT}
proved

\begin{lemma}\label{ht}
Let $\Sigma$ be a standard bitableau. Then
$$
\gamma_\delta(\Sigma)\neq 0 \implies  \gamma_\delta(\KRS(\Sigma))\neq 0
$$
\end{lemma}

\noindent and this implies that

\begin{theorem}\label{ht-2} The ideal $I(\delta,X)$ is G-KRS.
\end{theorem}

Note that from \ref{ht-2} one has that
$$
\gamma_\delta(\Sigma)\neq 0 \iff \gamma_\delta(\KRS(\Sigma))\neq 0.
$$

There are many natural questions concerning the function $\gamma_\delta$
and related ideals. For instance:

\begin{questions}
Let $J(\delta,k)=I_k(\gamma_\delta)$, that is,  the ideal generated by the
bitableaux $\Delta$ such that $\gamma_\delta(\Delta)\geq k$, and let
$B(\delta,k)$ be the set of the standard bita\-bleaux $\Sigma$ with
$\gamma_\delta(\Sigma)\geq k$.
\begin{itemize}
\item[(a)] Is $\gamma_\delta$ a KRS-invariant?
\item[(b)] Is $B(\delta,k)$ a basis of $J(\delta,k)$?  i.e.\ is
$\gamma_\delta$ str-monotone?
\item[(c)] Is $J(\delta,k)$ equal to $I(\delta,X)^{(k)}$?  The inclusion
$ J(\delta,k)\subseteq I(\delta,X)^{(k)}$ holds since the symbolic powers form
a filtration and for a single minor $\mu$ is not difficult to see that $\mu \in
I(\delta,X)^{(k)}$ where $k=\gamma_\delta(\mu)$.
\item[(d)] Is  $\ini(I(\delta,X)^{(k)})=\ini(I(\delta,X))^{(k)}$ ?
\end{itemize}
\end{questions}


\begin{thebibliography}{15.}

\bibitem{AK1} S.S. Abhyankar and D.M. Kulkarni. {\em Bijection between indexed
monomials and standard bitableaux}. Discrete Math. {\bf 79} (1990), 1--48.

\bibitem{AK2} S.S. Abhyankar and D.M. Kulkarni. {\em Coinsertion and
standard bitableaux}.) Discrete Math. {\bf 85} (1990), 115--166.

\bibitem{Ba} C. Bae\c tica. {\em Rees algebra of ideals generated by
pfaffians}. Commun. Algebra {\bf 26} (1998), 1769--1778.

\bibitem{BC} W. Bruns and A. Conca. {\em KRS
and powers of determinantal ideals}. Compositio Math. {\bf 111} (1998),
111--122.

\bibitem{BK} W. Bruns and M. Kwieci\'nski. {\em Generic graph construction
ideals and Greene's theorem}. Math. Z. {\bf 233} (2000), 115--126.

\bibitem{BV} W.~Bruns and U.~Vetter. {\em Determinantal rings}.
Lect.~Notes Math. {\bf 1327}, Springer 1988.

\bibitem{Co} A.~Conca. {\em Gr\"obner bases of ideals of minors of
a symmetric matrix}.  J.~Algebra {\bf 166} (1994), 406--421.

\bibitem{DeN} E. De
Negri. {\em ASL and Groebner bases theory for Pfaffians and monomial algebras}.
Dissertation, Universit\"at Essen (1996).

\bibitem{DEP} C.~De Concini, D.~Eisenbud and C.~Procesi. {\em Young diagrams
and determinantal varieties}. Invent. math. {\bf 56} (1980), 129--165.

\bibitem{F} W. Fulton. {\em Young tableaux}. Cambridge University Press 1997.

\bibitem{Gr} C. Greene. {\em An extension of Schensted's theorem}.
 Adv. Math. {\bf 14} (1974), 254--265.

\bibitem{HT} J.~Herzog and N.~V.~Trung. {\em Gr\"obner bases and
multiplicity of determinantal and pfaffian ideals}. Adv. in Math. {\bf 96}
(1992), 1--37.

\bibitem{K} D.E.~Knuth. {\em Permutations, matrices, and
generalized Young tableaux}. Pacific J. Math. {\bf 34} (1970), 709--727.

\bibitem{Sch} C.~Schensted. {\em Longest increasing and decreasing
subsequences}. Can.~J.~Math. {\bf 13} (1961), 179--191.

\bibitem{Sta} R.P. Stanley, {\em Enumerative Combinatorics}, Vol 2.
Cambridge University Press 1999.

\bibitem{Stu} B.~Sturmfels. {\em Gr\"obner bases and Stanley
decompositions of determinantal rings}. Math.~Z. {\bf 205} (1990), 137--144.

\bibitem{STV} B. Sturmfels, Ng\^o Viet Trung and W. Vogel.
{\em Bounds on degrees of projective schemes}. Math. Ann. {\bf 302} (1995),
417--432.

\end{thebibliography}
\end{document}